\documentclass[12pt,draft]{amsart}

\usepackage{amssymb,amscd,amsmath}
\usepackage[english]{babel}

\theoremstyle{plain}
\newtheorem{theorem}{Theorem}[section]
\newtheorem{lemma}[theorem]{Lemma}

\newtheorem{proposition}[theorem]{Proposition}

\theoremstyle{definition}

\newtheorem{example}[theorem]{Example}

\numberwithin{equation}{section}

\newcommand\NN{{\mathbb N}}

\newcommand\frk{\mathfrak}

\newcommand\mm{{\frk m}}

\newcommand\nn{{\frk n}}

\newcommand\opn[2]{%
    \newcommand{#1}{\operatorname{#2}}}

\opn\chara{char}
\opn\projdim{proj\,dim}
\opn\depth{depth}
\opn\rank{rank}
\opn\supp{supp}
\opn\Ass{Ass}
\opn\Proj{Proj}
\opn\Spec{Spec}
\opn\Soc{Soc}
\opn\Shad{Shad}
\opn\Span{span}
\opn\Hs{Hilb}
\opn\ini{in}
\opn\gin{Gin}
\opn\grade{grade}
\opn\Ker{Ker}
\opn\Coker{Coker}

\opn\im{Im}
\opn\Hom{Hom}
\opn\Tor{Tor}
\opn\Ext{Ext}
\opn\id{id}
\opn\Id{Id}
\opn\Homst{^*Hom}
\opn\Extst{^*Ext}
\opn\Gamst{^*\Gamma}
\opn\Hst{^*H}

\opn\lk{lk}
\opn\st{st}

\let\:=\colon
\let\phi=\varphi

\newcommand\thh{^{th}}

\newcommand\pnt{{\raise0.5mm\hbox{\large\bf.}}}

\newcommand\ra{\rightarrow}
\newcommand\Ra{\Rightarrow}
\newcommand\Lra{\Leftrightarrow}
\newcommand\p{^\prime}

\newcommand\ov{\overline}
\newcommand\lex{^{lex}}
\newcommand\sat{^{sat}}
\opn\po{pol}
\opn\pol{^{\bf p}} 

\begin{document}
\title {Ideals with maximal local cohomology modules}
\author{Enrico Sbarra}
\address{Enrico Sbarra\\ DSM - Universit\'a di Trieste\\
Via A. Valerio 12/1\\
I-34127 Trieste}

\medskip
\email{sbarra@dsm.univ.trieste.it}

\maketitle

\section*{Introduction and Notations}
This paper finds its motivation in the pursuit of ideals whose local cohomology modules have maximal
Hilbert functions. In \cite{S1}, \cite{S2} we proved that the lexicographic (resp. squarefree 
lexicographic) ideal of a family of graded (resp. squarefree) ideals with assigned Hilbert function
 provides sharp upper bounds 
for the local cohomology modules of any of the ideals of the family. Moreover these bounds are 
determined explicitely in terms of  the Hilbert function, which is the specified starting data. 
In the 
present paper a characterization of the class of ideals with the desired property is accomplished.\\
In order to be more precise we set some notation to be used henceforth.\\  
Let $R\doteq K[X_1,\ldots,X_n]$ denote the polynomial ring in $n$ variables over a field $K$ of 
characteristic $0$ with 
its standard grading, $\mm\doteq (X_1,\ldots,X_n)$ the maximal homogeneous ideal of $R$, 
$I\subset R$ a homogeneous ideal and $I\lex$ its lexicographic ideal. 

The canonical module of $R$ will be denoted by  $\omega_R\simeq R(-n)$. 
If $M$ stands for a  graded $R$-module, then $\Hs(M,t)$ will denote its Hilbert series in terms of $t$. 
The local cohomology modules $H^i_\mm(M)$ of $M$ will be considered with support on the maximal graded
ideal $\mm$ and with their natural grading. We write $h^i(M)_j$ for the dimension as a $K$-vector space
of $H^i_\mm(M)_j$. The dual of the local cohomology modules according to the Local Duality Theorem
will be denoted with $E^i(M)\doteq \Ext^i_R(M,\omega_R)$.
We write $I\sat\doteq I:\mm^\infty$ for the saturation of an ideal $I$ with respect to $\mm$.\\

A well known theorem of Bigatti-Hulett-Pardue states that in the family of graded ideals with a 
given Hilbert function the lexicographic ideal has the greatest Betti numbers. The class of ideals
with all of the Betti numbers equal to those (i.e. with the same resolution as that) of the 
lexicographic ideal 
are characterized
in \cite{AHH}. These ideals are the so-called {\it Gotzmann ideals}. We recall that an ideal is called
Gotzmann iff in each degree it has the same numbers of generators as its associated lexicographic ideal.
Note that the definition can be re-read as follows: Gotzmann ideals have the same $0\thh$ graded Betti
numbers as those of the associated lexicographic ideal. 
The result of \cite{AHH} shows then that the maximality of all
the other graded Betti numbers is forced by that of the $0\thh$ ones.\\
It is worth to point out that a similiar behaviour underlies our situation as well, where Gotzmann 
ideals will play again some important role. We state now the main result of this paper.

\begin{theorem}\label{main}
For any graded module $I$, the following are equivalent conditions:
\begin{itemize}
\item[(i)]
$(I\sat)\lex=(I\lex)\sat$;
\item[(ii)]
$h^i(R/I)_j=h^i(R/I\lex)_j$, for all $i,j$.
\end{itemize}
\end{theorem}
 
\indent The author would like to express his gratitude to Aldo Conca for many fruitful discussions
about the contents of this paper.

\section{Graded ideals}
Let us start by observing that there is an obvious inclusion in the equality $(i)$ of Theorem 
\ref{main}.
In fact one has that, for any ideal $I$, 
\begin{equation}\label{ink}
(I\sat)\lex\subseteq(I\lex)\sat.
\end{equation} 
One can prove this fact by an iterated use of $(I:\mm)\lex\subseteq I\lex:\mm$, which is easy.
Observe that \eqref{ink} provides an information about the Hilbert function of the $0\thh$ local 
cohomology modules. Since $H^0_\mm(R/I)\simeq I\sat/I$, 
$h^0(R/I)_j=\dim_KI\sat_j-\dim_KI_j=\dim_K(I\sat)_j\lex-\dim_KI\lex_j$, and by virtue of the above 
inclusion, this is less than or equal to $\dim_K(I\lex)\sat_j-\dim_KI\lex_j=h^0(R/I\lex)_j$. One sees
immediately that equality holds iff $(I\sat)\lex=(I\lex)\sat$.\\
In other words one could state Condition (i) as follows:
$$(i^*):\; h^0(R/I)_j=h^0(R/I\lex)_j, \hbox{ for all } j.$$ This also shows that $(ii)\Ra (i)$.
Thus, the theorem states that, as it happens in the context of  Betti numbers, 
the equality of the $0\thh$ local cohomology forces the equality of any other. One can wonder if this
sort of rigidity behaviour is to be expected more generally, i.e. is it true that if $h^i(R/I)_k=
h^i(R/I\lex)_k$ for some $i$ and all $k$, then $h^j(R/I)_k=h^j(R/I\lex)_k$ for all $j\geq i$ and all $k$? 
At this moment we have no evidence for answering the last question.


Let us denote by $\gin(I)$ the {\it generic initial ideal of $I$} with respect to the reverse 
lexicographical order induced by the assignment $X_1>X_2>\ldots >X_n$. 
In \cite{HS} we studied the problem of characterizing those ideals such that $h^i(R/I)_j=
h^i(R/\gin(I))_j$. For the reader's sake we recall the main theorem here:

\begin{theorem}\label{nostro}
Let $M$ be a finitely generated graded $R$-module with graded free presentation $M=F/U$. The following
conditions are equivalent.
\begin{itemize}
\item[(a)]
$F/U$ is sequentially CM;
\item[(b)]
for all $i\geq 0$ and all $j$ one has $h^i(F/U)_j=h^i(F/\gin(U))_j$.
\end{itemize}
\end{theorem}

In general it holds that $h^i(R/I)_j\leq h^i(R/\gin(I))_j\leq h^i(R/I\lex)_j$ (see \cite{S1}).
Thus, the class of ideals we are searching for must have the property that $R/I$ is sequentially
CM. Therefore one may state a third condition, wondering if this is equivalent to those of 
Theorem \ref{main}:
$$(iii):\; R/I \hbox{ is sequentially CM and } \gin(I)=I\lex.$$
By virtue of the above theorem one sees that $(iii)\Ra (ii)$. In fact, if $R/I$ is sequentially
CM, for all $i,j$ one has $h^i(R/I)_j=h^i(R/\gin(I))_j$, where the latter is equal to $h^i(R/I\lex)_j$,
since $\gin(I)=I\lex$.\\
On the other hand, one sees that in general $(iii)$ needs not to be implied by $(i)$. 
First let us prove an easy lemma.

\begin{lemma}\label{ginscambio}
Let $I$ be a homogeneous ideal and $\gin(I)$ its generic initial ideal. If  $(I\sat)\lex=
(I\lex)\sat$ then $(\gin(I)\sat)\lex=(\gin(I)\lex)\sat$.
\end{lemma}
\begin{proof}
As observed at the beginning of this section, one  inclusion is trivially true:
$(\gin(I)\sat)\lex\subseteq 
(\gin(I)\lex)\sat$. Thus, it is enough to show that the two ideals have the same Hilbert
function. Since $I$ and $\gin(I)$ have the same Hilbert function and therefore same lexicographic 
ideal,
one has $(\gin(I)\lex)\sat=(I\lex)\sat=(I\sat)\lex$, by hypothesis. Recall now that 
$\gin(I\sat)=\gin(I)\sat$. Therefore, $H(I\sat,t)=H(\gin(I\sat),t)=H(\gin(I)\sat,t)$, and as a 
consequence $H(I\sat,t)=H((\gin(I)\sat)\lex,t)$. Finally, we deduce that
$$\begin{array}{ll}
H((\gin(I)\lex)\sat,t)&=H((I\sat)\lex,t)=H(I\sat,t)\\
&=H((\gin(I)\sat)\lex,t),
\end{array}$$
 which yields the desired conclusion.
\end{proof}

By virtue of the above lemma, we need now an example of a strongly stable ideal, which is not 
lexicographic, but satisfies $(i)$. This is provided in what follows.

\begin{example}
Let $I=(x^2,xy,y^2,xz^2,yz^2)$ be an ideal of $K[x,y,z]$.\\
It is easy to verify that $I$ is strongly stable. The result of the computation of
the lex-ideal is $I\lex=(x^2, xy, xz, y^3, y^2z, yz^2)$. 
Thus, $I\neq I\lex$ and $(I\sat)\lex=(I\lex)\sat=(x,y)$.
\end{example}

Next, some lemmata which illustrate our hypothesis and characterize it. We shall prove that an ideal
with the exchange property is sequentially CM.

\begin{lemma}\label{uno}
Let $I$ be a homogeneous ideal. Then 
$$I=I\sat \hbox{ and } (I\sat)\lex=(I\lex)\sat \Lra I\lex=(I\lex)\sat.$$
\end{lemma}
\begin{proof}
``$\Ra$'': It is immediately seen.\\
``$\Leftarrow$'': Since $I\lex=(I\lex)\sat$, one has that $H^0_\mm(R/I\lex)=0$ and, by virtue
of \cite{S1}, Theorem 5.4, also $H^0_\mm(R/I)=0$, i.e. $I=I\sat$. Now, as already observed
elsewhere, it suffices to show that $H((I\sat)\lex,t)=H((I\lex)\sat,t)$. For this purpose we 
notice that $H((I\sat)\lex,t)$ is equal to $H(I\sat,t)=H(I,t)=H(I\lex,t)=H((I\lex)\sat),t)$, 
and we are done.
\end{proof}

\begin{lemma}\label{due}
Let $I$ be a homogeneous ideal. If $I=I\sat$ is Gotzmann then $(I\sat)\lex=(I\lex)\sat$.
\end{lemma}
\begin{proof}
By hypothesis  $\depth R/I>1$ and since  $I$ is a Gotzmann ideal, it has the same
resolution as $I\lex$. Therefore, also $R/I\lex$ has positive depth. 
This implies that $I\lex=(I\lex)\sat$ and the thesis follows immediately.
\end{proof}

The stronger counterpart of the above lemma is the following:
Let $I$ be a homogeneous ideal such that $(I\sat)\lex=(I\lex)\sat$. Then $I\sat$ is Gotzmann.
One can see that this last statement is equivalent to that of the following lemma.

\begin{lemma}\label{GotzMe}
Any ideal $I$ such that $I\lex$ is saturated is a Gotzmann ideal. 
\end{lemma}

Before we start to prove the latter fact, it is worth to 
underline that if the lexicographic ideal is saturated, the Hilbert function has a very rigid 
behaviour. In fact, saying that any ideal  associated with that lexicographic ideal is Gotzmann implies
that any of these ideals has the same resolution as the lex-ideal.\\
A saturated lexicographic ideal has indeed a very special structure and its generators can be 
described explicitely in terms of the Hilbert polynomial of $R/I$ by means of a vector $v$ of integers 
in a way that we are going to recall for the reader's sake.\\
Let $P_{R/I}(X)=\tbinom{X+a_1}{a_1}+\tbinom{X+a_2-1}{a_2}+\ldots+\tbinom{X+a_l-(l-1)}{a_l}$, with
$a_1\geq a_2\geq \ldots a_l\geq 0$ be a representation of the Hilbert polynomial,  also
referred to as its {\it Gotzmann representation}. Let now $v_i\doteq |\{j\: n-a_j-1=i\}|$, for 
$i=1,\ldots,n-1$ and observe that $v_i$ represents the exponent of the variable $X_i$ in the last 
monomial of highest degree of the minimal set of generators of $I\lex$. Let us also set 
$h$ to be the maximum index of a non-zero $v_i$. Then, the minimal set of generators of $(I\lex)\sat$
is the set 
$$\{X_1^{v_1+1},X_1^{v_1}X_2^{v_2+1},\ldots,X_1^{v_1}\cdot\ldots\cdot X_{h-1}^{v_{h-1}+1},
X_1^{v_1}\cdot\ldots\cdot X_{h-1}^{v_{h-1}}X_h^{v_h}\},$$
where, according to our settings, $v_i\geq 0$, for $i=1,\ldots,h-1$ and $v_h>0$. We also recall that
the vanishing of the local cohomology modules $H^i_\mm(R/I\lex)$ is determined by the vanishing 
of the $(n-i)\thh$ entry of the vector $v=(v_1,\ldots,v_h)$ (cf. \cite{S1}, Proposition 6.6).\\ 

\noindent
We may now continue our proof using an induction argument on $h$. If $h=1$, then $I\lex$ is simply
the principal ideal $(X_1^a)$, for some $a\in\NN_{>0}$, and there is nothing to prove. Suppose now
the thesis proven for any ideal such that the length of the vector $v$ is $h-1$. There are two possible
cases. If $v_1=0$, then the ideal $I\lex$  contains the linear form $X_1$, and consequently 
$I$ contains a linear form, let us say $l$. Thus $I\lex=(X_1,J)$ and $I=(l,I\p)$, for some ideals
$J=J\p R$, where $J\p$ is the saturated lexicographic ideal in $R\p\doteq K[X_2,\ldots,X_n]$ represented by the vector $(v_2,\ldots,v_h)$, and $I\p\subset R$.
Clearly $X_1$ is $R/J$-regular, and we also may assume that $l$ is $R/I\p$-regular. From this fact
one deduces that $J\p$ is the lexicographic ideal associated with $I\p R\p$, and one can use the inductive hypothesis to reach the conclusion.\\ 
Otherwise, if $v_1>0$, we observe that $\grade_{I\lex}(R)=1$ since $E^1(R/I\lex)$ $\neq 0$,
and that $I\lex=X_1^{v_1}J$, where $J$ is a saturated lexicographic ideal
represented by the vector $(0,v_2,\ldots,v_h)$ in $K[X_1,\ldots,X_n]$. 
Moreover $v_1$ equals the multiplicity $e$ of $R/I\lex$,
which is the same as that of $R/I$. We also have that $\grade_I(R)=1$, since the dimension of $R/I$ is the same as that of $R/I\lex$, therefore $I$ can be written as a product of a polynomial $f$
times an ideal $I\p$. Observe that $\deg f$ equals the multiplicity of $R/I$, which is $v_1$. For 
showing this simple fact, let us look at the short exact sequence $0 \ra fR/I \ra R/I \ra R/fR \ra 0$.
It is easy to see that the multiplicity of $R/fR$ is $\deg f$, since the $h$-vector of $R/fR$ is
$\sum_{i=1}^{\deg f}t^i$, while its dimension is $n-1$. On the other hand the Hilbert function
of the left-hand side module is up to a shift that of the module $R/J$, whose dimension is less than 
$n-1$ and therefore cannot contribute to the multiplicity of $R/I$.\\
Since $t^{\deg f}\Hs(I\p,t)=\Hs(I,t)=\Hs(I\lex,t)=t^{v_1}\Hs(J,t)$, we deduce that $J$ is the 
lexicographic ideal associated with $I\p$, and the proof of the lemma is complete by the use of the 
previous case. 

\begin{proposition}\label{sCM2}
Let $I$ be an ideal such that $(I\sat)\lex=(I\lex)\sat$. Then $R/I\sat$ is sequentially CM.
\end{proposition}

It is now convenient to  recall the definition of  sequentially Cohen-Macaulay modules.
A finitely generated graded $R$-module $M$ is said to be {\it sequentially Cohen-Macaulay} if there
exists a finite filtration $0=M_0\subset M_1 \subset M_2 \subset \ldots\subset M_r=M$ of $M$ by graded
submodules of $M$ such that each quotient $M_i/M_{i-1}$ is CM, and $\dim (M_1/M_0)<\dim (M_2/M_1)
<\ldots \dim (M_r/M_{r-1}).$\\
Sequentially CM modules have an interesting characterization in terms
of their homology, as illustrated by a theorem of Peskine (cf. \cite{HS}, Theorem 1.4).
Peskine's result asserts that a module $M$ is sequentially CM iff for all $0\leq i\leq \dim M$ the
modules $E^{n-i}(M)$ are either $0$ or CM of dimension $i$. 
For a more complete overview of the properties of sequentially CM modules we refer to \cite{HS}, and
we proceed by proving Proposition \ref{sCM2}.

\begin{proof}
Let us assume that $I\lex$ is saturated and let us prove that $R/I$ is sequentially CM. 
We shall use the same notation as above. In particular the vector that
determines $I\lex$ will be denoted by $v$, and $v_i$, for $i=1,\ldots,h$,
will denote its entries.\\
The idea is the same as that of the proof of the above lemma, by induction on the index $h$. If $h=1$, 
then $I$ and $I\lex$ are principal, therefore they are Cohen-Macaulay, and Cohen-Macaulay modules
are sequentially CM.\\
If $h>1$ we may assume without loss of generality that $v_1$ is $0$ and therefore that $I$ and $I\lex$
contain the linear form $X_1$. The following is an application of the graded
version of the Rees' Lemma. Let $I=(X_1,I\p)$ and $I\lex=(X_1,J)$, where
$J=J\p R$ and $J\p$ is the (saturated) lexicographic ideal associated with
$I\p(R/(X_1))$ determined by the vector $(v_2,\ldots,v_h)$ 
in $K[X_2,\ldots,X_n]$. 
Thus we have graded isomorphisms for all $i\geq 1$
$$\begin{array}{ll}
\Ext^i_R(R/I,R)&\simeq \Ext^i_R(R/(X_1,I\p),R)\\
&\simeq \Ext_{R/(X_1)}^{i-1}((R/(X_1))/(I\p(R/(X_1))),R/(X_1))(1),\\
\end{array}$$ which, by induction, 
is either $0$ or Cohen-Macaulay of dimension $(n-1)-(i-1)=n-i$. By virtue
of the aforementioned homological characterization of sequentially 
CM modules this is equivalent to the thesis.
\end{proof}

As an immediate consequence of the above observations, we deduce the property we were interested in.

\begin{proposition}\label{sCM}
Let $I$ be an ideal such that $(I\sat)\lex=(I\lex)\sat$. Then $R/I$ is sequentially CM.
\end{proposition}
\begin{proof}
By virtue of the above proposition it is enough to notice that an $R$-module $M$ is sequentially
CM iff $M/H^0_\mm(M)$ is sequentially CM. In our case $M=R/I$ and $M/H^0_\mm(M)\simeq (R/I)/(I\sat/I)\simeq R/I\sat$.
\end{proof}

Observe that if $I$ is a proper cyclic ideal, its lexicographic ideal is $(X_1^d)$ for some positive integer $d$, which is saturated. Moreover, $R/I$ is
CM of dimension $n-1$, and its only non-vanishing $\Ext$-group is the first one. 
This is isomorphic to a shifted copy of $R/I$ itself, and therefore cyclic.
Applying this observation to the inductive argument of the proof of Proposition
\ref{sCM2}, one deduces the following.

\begin{proposition}
Let $I$ be an ideal such that 
$(I\sat)\lex=(I\lex)\sat$, then all of its $\Ext$-groups except possibly 
the $n\thh$-one are cyclic.
\end{proposition}

We now prove the main result for strongly stable ideals. Recall that a monomial ideal $I$ is said to be 
{\it strongly stable} if, for every $u\in I$, one has $X_iu/X_j\in I$ for all $X_j|u$ and $i<j$.
For a strongly stable ideal one has that $I\sat=I:(X_n)^\infty$. In particular, one knows that
$R/I$ has positive depth iff $X_n$ does not appear in any of the monomials which minimally generate $I$.
Suppose now that $I$ is strongly stable and that $X_n$ is $R/I\lex$-regular. Then, if we denote
by $\ov{\cdot}$ the equivalence classes in the quotient of the polynomial ring by the last 
variable, one has $\ov{I\lex}=\ov{I}\lex$. Let us give a quick explanation for this fact. Since $I\lex$
is strongly stable as well and $X_n$ is $R/I\lex$-regular, none of the monomial of the minimal set of
generators of $I\lex$ contain the last variable, thus $\ov{I\lex}$ is  a lexicographic ideal in the
variables $X_1,\ldots,X_{n-1}$. Since $\depth R/I\lex\leq \depth R/I$, and $I$ is a strongly stable 
ideal, for the  reason explained above  $X_n$ is also $R/I$-regular. Therefore we can control the 
behaviour of the Hilbert function when passing to the quotient by the last variable, and the 
conclusion follows easily.

The following is a technical lemma about local cohomology and we refer to \cite{BS} or \cite{BH}
for more details about the Local Duality Theorem and the properties of the canonical module.

\begin{lemma}\label{lemmalemma}
Let $I\subset R$ an ideal of $R$ and let $S\doteq R[X]$, with maximal graded ideal $\nn$. 
The following  graded isomorphism of $R$-modules holds.
$$H^i_\nn(S/IS)\simeq \Hom_R(S,H^{i-1}_\mm(R/I))(1).$$
\end{lemma}
\begin{proof}
The relation $\omega_S=(\omega_R\otimes_R S)(-1)$ between the canonical modules of $S$ and $R$ is
well known. By the Local Duality Theorem, one has that $H^i_\nn(S/IS)$ is the dual of 
$\Ext_S^{n+1-i}(S/IS,\omega_S)$, which is defined to be $\Hom_K(\Ext_S^{n+1-i}(S/IS,\omega_S),K)$.
After writing $S=S\otimes_RR$ and substituting $\omega_S$ by means of the formula written above, 
one obtains that the latter is isomorphic to
$\Hom_K(\Ext_{S\otimes_RR}^{n+1-i}(R/I\otimes_RS,(\omega_R\otimes_RS)(-1)),K)\simeq
 \Hom_K(\Ext_{S\otimes_RR}^{n+1-i}(R/I\otimes_RS,\omega_R\otimes_RS),K)(1)$ and, since $S$ is $R$-flat,
this is isomorphic to $\Hom_K(S\otimes_R\Ext_R^{n+1-i}(R/I,\omega_R),K)(1)$. Using a well-known
formula in homological algebra and applying the Local Duality Theorem again, one finally deduces
$$\begin{array}{ll}
H^i_\nn(S/IS)&\simeq \Hom_R(S,\Hom_K(\Ext_R^{n+1-i}(R/I,\omega_R),K))(1)\\
&\simeq\Hom_R(S,H_\mm^{i-1}(R/I))(1),
\end{array}$$ as required.
\end{proof}

\begin{lemma}\label{induttivo}
Let $S=R[X]$ a polynomial ring in one indeterminate over $R$ with graded maximal ideal $\nn$. 
Let $I$ be an ideal of $R$. Then, for all $i,j$,  
$$\dim_KH^i_\nn(S/IS)_j=\sum_{h\geq j}\dim_KH^{i-1}_\mm(R/I)_{h+1}.$$
\end{lemma}
\begin{proof}
The proof is based on the above lemma and on some considerations about the $S$-module structure 
of $\Hom_R(S,N)$, where $N$ is an arbitrary $R$-module. Let us consider the $R$-linear application
$\cdot:\;S\times \Hom_R(S,N)\ra \Hom_R(S,N)$ defined by $\cdot(s,\phi)(t)\doteq s\cdot\phi(t)\doteq
\phi(st)$. This map  endows $\Hom_R(S,N)$ with an $S$-module structure.\\

\noindent
Let now $N$ and $S$ be  graded. One defines a graded structure of $S$-module on $\Hom_R(S,N)$
as follows. We set $$\Hom_R(S,N)_i\doteq\{\phi\in\Hom_R(S,N)| \phi(S_j)\subset N_{i+j}, \hbox{for all 
} j\}$$ to be the $i\thh$ graded component of $\Hom_R(S,N)$.\\
Observe that, if $\phi\in\Hom_R(S,N)_i\cap\Hom_R(S,N)_j$ and $i\neq j$, then $\phi(S_k)\subset
N_{i+k}\cap N_{j+k}=(0)$, i.e. $\phi=0$. If $\phi\in\Hom_R(S,N)_i$, $s$ is an element of $S_j$ and 
$t\in S_k$ then $s\phi(t)=\phi(st)\in N_{i+(k+j)}$, i.e. $s\phi\in \Hom_R(S,N)_{i+j}$. Finally, 
one can verify that  $\Hom_R(S,N)\subset \oplus_i \Hom_R(S,N)_i$.\\

Let now $S$ be as in the hypothesis and consider the homogeneous isomorphism of graded $S$-modules
$\alpha$
$$\Hom_R(S,N)\stackrel{\alpha}{\ra}\Pi_{j\geq 0}Nx^{-j}$$
that maps an element $\phi$ of $\Hom_R(S,N)$ into $(\ldots,\phi(x^j)x^{-j},\ldots)$.\\
Thus, $\Hom_R(S,N)_i\simeq (\Pi_{j\geq 0}Nx^{-j})_i\simeq \oplus_{h\geq i}N_h$, and if $N$ is artinian
one can deduce that the dimension as a $K$-vector space of $\Hom_S(R,N)$ is just the sum of 
the dimensions of the graded components $N_h$ of $N$ with $h\geq i$. Now the conclusion follows 
from  Lemma \ref{lemmalemma}. 
\end{proof}

\begin{proposition}\label{strongly}
Let $I$ be a strongly stable ideal such that $(I\sat)\lex=(I\lex)\sat$. Then
$$h^i(R/I)_j=h^i(R/I\lex)_j, \hbox{ for all } i,j.$$
\end{proposition}
\begin{proof}
As we noticed already several times, the hypothesis is equivalent to saying that $h^0(R/I)_j=
h^0(R/I\lex)_j$ for all $j$. We may then assume that $I$ is saturated, i.e. that $\depth R/I$ is
positive. By induction on $n$ we also suppose  the thesis to be true  for any 
strongly stable ideal on a polynomial ring with less than $n$ variables. Bearing in mind the remarks 
before Lemma \ref{induttivo}, the conclusion is a straightforward application of the latter.   
\end{proof}

\begin{proof}[Proof of Theorem \ref{main}, $(i) \Ra (ii)$]
By Proposition \ref{sCM}, we know that $R/I$ is sequentially CM. If $\gin(I)=I\lex$, we achieve the
conclusion immediately for what we said before Lemma \ref{ginscambio}. Otherwise, by virtue of
Lemma \ref{ginscambio} we may assume without loss of generality that $I$ is strongly stable, by 
substituting it with its generic initial, since $h^i(R/I)_j=h^i(R/\gin(I)_j$. The thesis is thus an
application of Proposition \ref{strongly}.
\end{proof}





\begin{thebibliography}{1111}

\bibitem {AHH}
A. Aramova, J. Herzog, T. Hibi. Ideals with stable Betti numbers. {\em Adv. Math.} {\bf 152} (1) (2000),
72-77.

\bibitem {Bi}
A. Bigatti. Upper bounds for the Betti numbers of a given Hilbert function. {\em Communications in
 Algebra} {\bf 21}(7)
(1993), 2317-2334.

\bibitem {BS}
M. Brodmann, R. Sharp. {\it Local Cohomology}. Cambridge University Press, Cambridge, 1998.
 
\bibitem {BH}
W. Bruns, J. Herzog. {\it Cohen-Macaulay rings}. Cambridge University Press, Cambridge, 1998.

\bibitem {Hu}
H. Hulett. Maximum Betti numbers of homogeneous ideals with a given Hilbert function. {\em Communications in Algebra} 
{\bf 21}(7) (1993), 2335-2350.

\bibitem {HS}
J. Herzog, E. Sbarra. Sequentially Cohen-Macaulay modules and local cohomology. {\it To appear in}
Proceedings of the Conference of the International Colloquium on Algebra, Arithmetic and Geometry,
TIFR, Mumbai, January 4-12, 2000.

\bibitem {P}
K. Pardue. Deformation classes of graded modules and maximal Betti numbers. {\em Illinois Journal of Mathematics}
{\bf 40}(4) (Winter 1996), 564-585.

\bibitem {S1}
E. Sbarra. Upper Bounds for local cohomology for rings with given Hilbert function. {\em Communications
in Algebra} {\bf 29}(12) (2001) 5383-5409. 

\bibitem {S2}
E. Sbarra. On the structure of $\Ext$ groups of strongly stable ideals. {\em Geometric and combinatorial
aspects of commutative algebra}, Lect. Notes in Pure and Applied Mathematics, {\bf 217}, 345-352,
Dekker 2001.
\end{thebibliography}
\end{document}